\documentclass[a4paper,12pt]{article}
\usepackage{amssymb,amsmath}

\def\endproof{\hfill$\square$}

\newlength{\argwidth}

\makeatletter
\def\overline{\@ifnextchar [\myoverline{\@@overline}}
\def\myoverline[#1]#2{\mathchoice{\xoverline{#1}{#2}\displaystyle}%
{\xoverline{#1}{#2}\textstyle}{\xoverline{#1}{#2}\scriptstyle}%
{\xoverline{#1}{#2}\scriptscriptstyle}}

\def\xoverline#1#2#3{\settowidth{\argwidth}{$#3#2$}%
\lefteqn{\hbox to\argwidth{\hfil\yoverline{#1}{#2}{#3}\hfil}}#2}

\def\yoverline#1#2#3{$\@@overline{\vphantom{#3#2}\hspace{#1\argwidth}}$}
\makeatother

\def\bar#1{\overline[0.7]{#1}}

\makeatletter

\renewcommand{\@begintheorem}[2]{\begin{trivlist}\it\item[\hspace{\labelsep}{\bf #1\ #2.}]}
\renewcommand{\@opargbegintheorem}[3]{\begin{trivlist}\it\item[\hspace{\labelsep}{\bf #1\ #2\ (#3).}]}
\renewcommand{\@endtheorem}{\end{trivlist}}

\def\section{\@ifstar{\sectione}{\sections}}
\def\sectione#1{\begin{center} \bigskip \bf \large #1 \end{center}}
\def\sections#1{\refstepcounter{section}%
\begin{center} \bigskip \bf \large \thesection. #1 \end{center}}

\makeatother

\def\ep{\varepsilon}

\def\al{\alpha}

\def\sg{\sigma}

\def\lm{\lambda}
\def\Lm{\Lambda}
\def\om{\omega}

\def\N{\mathbb N}
\def\vp{\varphi}

\def\cc{$\cdots$}
\def\<{\left<}
\def\>{\right>}
\def\N{\mathbb N}

\def\st{\subset}
\def\tleq{\trianglelefteq}
\def\treq{\trianglerighteq}

\def\upa{{\uparrow}}
\def\doa{{\downarrow}}
\def\stn{\varsubsetneq}

\newcommand{\diag}{\mathop{\rm diag} \nolimits}
\newcommand{\ch}{\mathop{\rm char} \nolimits}

\newcommand{\End}{\mathop{\rm End}}
\newcommand{\Hom}{\mathop{\rm Hom}}

\newcommand{\Ker}{\mathop{\rm Ker}}
\newcommand{\supp}{\mathop{\rm supp}}
\newcommand{\dual}{\mathop{\rm d\vphantom{u}}}
\newcommand{\im}{\mathop{\rm Im}}

\newtheorem{teo}{Theorem}
\newtheorem{lemma}{Lemma}
\newtheorem{sled}{Corollary}
\newtheorem{usl}{Condition}
\newtheorem{problem}{Problem}

\newenvironment{explicit}[1]{\par\medskip\noindent{\bf#1.} \it}{\medskip}

\def\sgn{\mathop{\rm sgn}}

\def\^#1{{}^{#1}\hspace{-0.4em}}
\def\|{\linebreak[0]}

\begin{document}

\begin{center}
\Large Finite basis property for some classes of irreducible
representations of the symmetric groups
\normalsize
\vskip 1.5em%
{V.\,V. Shchigolev}\\
{\bf e-mail} vkshch\char64 vens.ru
\end{center}

\makeatletter
\def\footnoterule{\kern-3\p@
  \hrule width 3truein \kern 2.6\p@}
\makeatother

\insert\footins{
}

\begin{abstract}
In this paper we study the possibility to define
irreducible representations of the symmetric groups
with the help of finitely many relations.
The existence of finite bases is established for
the classes of representations corresponding to two-part partitions and
to partitions from the fundamental alcove.
\end{abstract}

\section{Introduction}\label{s:0}

In this paper we continue the research begun in~\cite{Shchigolev9}.
Let us recall the basic notation used there.
Let $X_n=\{x_t^s:t\in\N,s=1$, \ldots, $n\}$ be a set of commuting
variables, $K$ be a field and $F_n=K[X_n]$
be the free commutative algebra with $1$.
Denote by $F_n^r$ the subspace of $F_n$
spanned by monomials of the form $x_{i_1}^{j_1}$\cc$x_{i_r}^{j_r}$,
where $\{i_1$, \ldots, $i_r\}=\{1$, \ldots, $r\}$.
The symmetric group $G(r)$ of degree $r$ acts on $F_n^r$ on the left
by the formula
$\pi x_{i_1}^{j_1}$\cc$x_{i_r}^{j_r}=x_{\pi i_1}^{j_1}$\cc$x_{\pi i_r}^{j_r}$.

Put $[x_{t_1}$, \ldots, $x_{t_k}]=%
\sum_{\sg\in G(k)}\sgn(\sg)x_{t_1}^{\sg(1)}$\cc$x_{t_k}^{\sg(k)}$,
where $k=1$, \ldots, $n$.
Let $\lm$ be a partition of $r$ into no more than $n$ parts.
Denote by $S^{\lm}$ the subspace of $F_n^r$
spanned by polynomials of the form
$\prod_{i=1}^{\lm_1}[x_{t_{i,1}}$, \ldots, $x_{t_{i,\mu_i}}]$,
where $\mu$ is the partition conjugate to $\lm$ and
$\{t_{i,j}:1\le i\le\lm_1,1\le j\le\mu_i\}=\{1$, \ldots, $r\}$.
For any positive integer $i$, denote by $\lm+(i^n)$
the partition $(\lm_1+i$, \ldots, $\lm_n+i)$ (see Section~\ref{s:2}).

Let $U$ be a $KG(r)$-submodule of $F_n^r$.
Denote by $U\upa$, or more precisely by $U\upa_n$ if we want to
underline the role of $n$, the $KG(r+n)$-submodule of $F_n^{r+n}$
generated by the subspace $U[x_{r+1}$, \ldots, $x_{r+n}]$.
Note that $U\st S^\lm$ implies $U\upa \st S^{\lm+(1^n)}$.

Conversely, let $V$ be a $KG(r+n)$-submodule of $F_n^{r+n}$.
Denote by $V\doa$ (or by $V\doa_n$) the set of all
polynomials $f\in F_n^r$ such that $f[x_{r+1}$, \ldots, $x_{r+n}]\in V$.
Clearly, $V\doa$ is a $KG(r)$-module.
Applying~\cite[Corollary~17.18]{James1}
one can show that $V\st S^{\lm+(1^n)}$ implies
$V\doa\st S^\lm$.

Let $U\upa^t$ or $U\doa^t$ denote the result of the
$t$-fold application of $\upa$ or $\doa$ respectively
to $U$, of course, if this is possible
(for $U\doa^t$).
Let us remark the following simple properties:
$U\doa\upa\st U\st U\upa\doa$ and $U\doa\upa\doa=U\doa$
for any submodule $U$ of $F_n^r$, where $r\ge n$.

In~\cite{Shchigolev9} we formulated the following problem:
\setcounter{problem}{1}\label{probl:1}
\begin{problem}
Let
$V_i\st S^{\lm+(q_i^n)}$, $i\in\N$,
where $0\le q_1<q_2<$\cc, be a sequence of submodules such that
$V_i\upa^{q_{i+1}-q_i}\st V_{i+1}$ for every $i\in\N$.
Does there exist any $N$ such that, for every $i>N$, we have
$V_N\upa^{q_i-q_N}=V_i$?
\end{problem}

We considered this problem as some special, though important, case of
the following problem:

\begin{explicit}{Problem~1}
Suppose $\ch{K}=p>0$.
Is every sequence of multilinear polynomials from $A$ finitely based?
\end{explicit}

Here $A$ is the free commutative $K$-algebra
with free generators $f_i(x_{t_1}$, \ldots, $x_{t_{n(i)}})$,
where $i=1$, \ldots, $k$ and $t_1$, \ldots, $t_{n(i)}\in\N$,
a polynomial of $A$ is considered to be multilinear
if it is multilinear with respect to the multiplicities of
the variables $x_j$ (for example, $f(x_1,x_2)f(x_2,x_3)$ is not a
multilinear polynomial) and implication is considered with respect to
$K$-linear actions, multiplication by elements of $A$
and renaming the variables $x_j$.

Problem~1 can be regarded as the analog of Specht's problem for forms.
In that case form relations play the role of identities.
An example of these relations is the well known Garnir relations.
In the case $\ch K=0$ the latter define the Specht modules and
thus the irreducible modules.
If $\ch K=p>0$ then this is not the case and
a problem about the existence of finite sets of relations that
define some classes of irreducible modules arises.
By the results of~\cite{Shchigolev9} it is natural to consider the
 classes of the form $D^{\lm+(i^n)}$, $i\in\N\cup\{0\}$,
where $\lm+(1^n)$ is a $p$-regular partition.
This choice corresponds to solving Problem~2 in the case where
$q_i=i-1$ and $V_i$ is the maximal submodule of $S^{\lm+({i-1}^n)}$
for any $i\in\N$.

The positive solution to the last problem for partitions $\lm$ that
satisfy Condition~\ref{usl:1} formulated in Section~\ref{s:2}
at page~\pageref{usl:1} is given by Theorem~\ref{t:1}.
According to the results of~\cite{Kleshchev1} and~\cite{Kleshchev2}
there are some partitions $\lm$ not satisfying Condition~\ref{usl:1}.
However, we know that it is satisfied for partitions $\lm$
containing no more than $n$ nonzero parts such that
$\lm_1-\lm_n\le p-n$ and for two-part partitions.
This is proved in Section~\ref{s:3} with the help of
the results of the papers~\cite{Mathieu1} and~\cite{Erdmann1}.
Note that we consider the case $n=2$ only for illustration,
since~\cite[Theorem~6]{Shchigolev9}
solves the partial case $n=2$ of Problem~\ref{probl:1}
without any restriction on the submodules $V_i$.

\section{Preliminaries}\label{s:1}

Let us fix an algebraically closed field $K$ of positive characteristic $p$
and an integer $n\ge2$.
Let $E$ be the linear $K$-space with basis $v_1$,\ldots,$v_n$
and $GL_n(K)$ be the group of all $K$-linear automorphisms of $E$.
The group $GL_n(K)$ is naturally endowed with the structure of
an algebraic $K$-group and we may consider polynomial (rational)
$GL_n(K)$-modules with respect to this structure.
Denote by $\Lm$ the set of all sequences of whole numbers of length $n$.
We write $\lm\treq\mu$, where $\lm,\mu\in\Lm$, if
we have $\lm_1+$\cc$+\lm_i\ge\mu_1+$\cc$+\mu_i$ for any $i=1$,\ldots,$n$.
For a polynomial $GL_n(K)$-module $V$ and a sequence of whole numbers
$\lm=(\lm_1$,\ldots,$\lm_n)$, we denote by $V^\lm$ the set of all
$v\in V$ such that $\diag(t_1$,\ldots,$t_n)v=t_1^{\lm_1}$\cc$t_n^{\lm_n}v$
for any $t_1$,\ldots,$t_n\in K^*$.
Clearly, $V^\lm$ is a $K$-linear space.
This space is called the {\it $\lm$-weight space} of $V$.
If $V^\lm\ne0$ then we call $\lm$ a {\it weight} of $V$.
Denote by $\Lm^+$ the subset of $\Lm$ consisting of
sequences $(\lm_1$,\ldots,$\lm_n)$ such that $\lm_1\ge$\cc$\ge\lm_n$.

In the papers~\cite{Donkin1} and~\cite{Ringel1} the theory of
the so-called tilting modules is developed.
Here we assume that any tilting module is a finite dimensional
polynomial $GL_n(K)$-module.
We will not give definitions of these objects.
Instead, let us enumerate their basic properties
necessary for the sequel.
\begin{enumerate}
\item\label{p:1}  $E$ is a tilting module
(by the definitions of~\cite{Donkin1} and~\cite{Ringel1}).
\item\label{p:2} The tensor product of tilting modules is
also a tilting module
(the first full proof in the general case was obtained
by O.\,Mathieu~\cite{Mathieu2}).
\item\label{p:3} For any $\lm\in\Lm^+$ there exists
an indecomposable module $P(\lm)$ such that $\dim P(\lm)^\lm=1$
and if $\dim P(\lm)^\mu>0$ for some $\mu\in\Lm$, then $\mu\tleq\lm$
(see~\cite[Theorem~(1.1)]{Donkin1}).
\item\label{p:4}
Any tilting module is a direct sum of tilting modules
isomorphic to $P(\lm)$, $\lm\in\Lm^+$
(see~\cite[Theorem~(1.1)]{Donkin1}).
\end{enumerate}

For any positive integer $r$, we denote by $G(r)$ the symmetric group of
degree $r$, i.e. the set of bijections $\sg$ of the set of positive integers
such that $\sg(i)=i$ for any $i>r$.
This realization of the symmetric groups implies the
inclusions $G(r)\st G(r')$ for $r\le r'$.
Let $\sgn(\sg)$ denote the sign of a permutation $\sg$.
Let $\Lm^+_r$ denote the subset of $\Lm^+$ consisting of the sequences
of nonnegative integers with sum $r$.
For any set $S$ we denote by $i_S$ the identity map from $S$ to itself.

Consider an arbitrary decomposition
$E^{\otimes r}=\bigoplus_{i=1}^{h_r}M_i(r)$
into a direct sum of indecomposable $GL_n(K)$-modules.
By properties~\ref{p:1} and~\ref{p:2} the module $E^{\otimes r}$ is tilting.
Therefore by property~\ref{p:4} and
the Krull-Schmidt theorem
the modules $M_i(r)$ are tilting indecomposable.
Put $A(r)=\End_{GL_n(K)}E^{\otimes r}$.
We have the natural projections $\pi_i^r:E^{\otimes r}\to M_i(r)$ and
embeddings $\ep_i^r:M_i(r)\to E^{\otimes r}$.
Define $e_i^r=\ep_i^r\pi_i^r$.
The endomorphisms $e_i^r$ are primitive idempotents
and we have $e_i^re_j^r=0$ for $i\neq j$.

One can see directly from the definition that any weight of $E^{\otimes r}$
has nonnegative entries with sum $r$.
Therefore we have the partition
$\{1$,\ldots,$h_r\}=\bigsqcup_{\lm\in\Lm^+_r}I_\lm$,
where $I_{\lm}$ are possibly empty sets and $M_i(r)\cong P(\lm)$
for $i\in I_\lm$.

For any pair $i,j\in I_\lm$, we choose some isomorphism $\tau_{i,j}^r$
from $M_j(r)$ to $M_i(r)$ in such a way that
$\tau_{i,i}^r=i_{M_i(r)}$ and $\tau_{i,j}^r\tau_{j,k}^r=\tau_{i,k}^r$
for any $i,j,k\in I_\lm$.

Let us consider the following idempotents $e_\lm=\sum_{i\in I_\lm}e_i^r$.
In addition, we put $e_\lm=0$ if $I_\lm$ is empty.
We have the representation $A(r)=S(r)\oplus J(r)$,
where $S(r)$ is the set of endomorphisms of the form
$\sum_{\lm\in\Lm^+_r}\sum_{i,j\in I_\lm}\al_{i,j}$
$\ep_i^r\tau_{i,j}^r\pi_j^r$, where $\al_{i,j}\in K$,
and $J(r)$ is the set of endomorphism of the form
$\sum_{i,j=1}^{h_r}\ep_i^r \vp_{i,j}\pi_j^r$, where
$\vp_{i,j}\in\Hom_{GL_n(K)}(M_j(r),M_i(r))$ are such that
any $\vp_{i,j}\tau_{j,i}^r$ is nilpotent if
$i,j \in I_\lm$ and $\lm\in\Lm^+_r$.
This decomposition easily follows from the Fitting lemma and the fact that
$K$ is algebraically closed and all modules under consideration are
finite dimensional.
We have
$$
A(r)=\bigoplus_{i,j=1}^{h_r}e_i^rA(r)e_j^r=
\bigoplus_{\lm\in\Lm^+_r}S(r)e_\lm\oplus\bigoplus_{\lm,\mu \in \Lm^+_r}
e_\lm J(r)e_\mu.
$$
In addition $S(r)$ is a subalgebra of $A(r)$
and $S(r)e_\lm=e_\lm S(r)=e_\lm S(r)e_\lm$ for any $\lm\in\Lm^+_r$.

Given $\pi\in G(r)$, let us consider the element $\sg_r(\pi)$ of $A(r)$
acting by the formula
$\sg_r(\pi)(v_{i_1}\otimes$ \cc $\otimes v_{i_r})=
v_{i_{\pi^{-1}(1)}}\otimes$ \cc $\otimes v_{i_{\pi^{-1}(r)}}$.
Extending $\sg_r$ to the whole of $KG(r)$ by linearity,
we obtain a homomorphism from $KG(r)$ to $A(r)$.
With the help of this homomorphism we will consider any
$A$-module as a $KG(r)$-module.

\begin{lemma}[De~Concini, Procesi~{\rm\cite{DeConcini_Procesi1}}]%
\label{l:0.25}
The map $\sg_r$ is an epimorphism.
If $r\le n$ then $\sg_r$ is an isomorphism
and if $r\ge n+1$ then $\Ker\sg_r$ is generated as a two-sided ideal
of $KG(r)$ by $\sum_{\sg\in KG(n+1)}\sgn(\sg)\sg$.
\end{lemma}

Let us call a sequence of integers $\lm=(\lm_1$, \ldots, $\lm_s)$ such that
$\lm_1\ge$\cc$\ge\lm_s\ge0$ and $\lm_1+$\cc$+\lm_s=r$
a {\it partition} of $r$.
We call this sequence {\it $p$-regular} if there are no
$p$ numbers $\lm_i$ equal to a positive integer.
Otherwise this sequence is called {\it $p$-singular}.
Any $p$-regular partition $\lm$ of $r$ labels
in the standard way (see~\cite{James1}) some
irreducible $KG(r)$-module denoted by $D^\lm$.
In addition, any irreducible $KG(r)$-module is isomorphic to some $D^\lm$.
Using the epimorphism $\sg_r$ defined above, one can make
any irreducible $KG(r)$-module $N$ isomorphic to some $D^\lm$ for
a $p$-regular $\lm\in\Lm^+_r$ into an $A(r)$-module by the formula
$am=xm$, where $m\in N$, $a\in A(r)$, $x\in KG(r)$ and $\sg_r(x)=a$.
This definition is correct, since by Lemma~\ref{l:0.25}
and the direct construction of $D^\lm$ we have $(\Ker\sg_r)D^\lm=0$.
By Lemma~\ref{l:0.25} one can easily check that any irreducible
$A(r)$-module is isomorphic to $D^\lm$ for some
$p$-regular $\lm\in\Lm^+_r$.

We say an irreducible representation $D^\mu$ corresponds to
an idempotent $e_\lm$ if $e_\lm D^\mu\ne0$.
It is obvious that $e_\lm m=m$ for any $m\in D^\mu$ and $e_\nu D^\mu=0$
for $\nu\ne\lm$.
From the general theory it follows that for any nonzero primitive idempotent
there is one up to isomorphism irreducible representation corresponding to
this idempotent.

\begin{lemma}\label{l:0.5}
The representation $D^\lm$, where $\lm$ is a $p$-regular partition
of $\Lm^+_r$, corresponds to $e_\lm$.
If $\lm$ is a $p$-singular partition, then $e_\lm=0$.
\end{lemma}
{\bf Proof.}
For brevity, we put
$G=G(r)$, $A=A(r)$, $M_i=M_i(r)$ and $e_i=e_i^r$.
Denote by $X$ the set consisting of $\lm\in\Lm^+_r$ such that $e_\lm\ne0$.
For $\lm\in X$ we denote by $t(\lm)$ the $p$-regular partition
of $\Lm^+_r$ such that $e_\lm D^{t(\lm)}\ne0$.
Clearly, $t$ is a bijection from $X$ to the subset of $\Lm^+_r$
consisting of $p$-regular partitions.

For any $\lm\in\Lm^+_r$ we put
$$
v(\lm)=c_1\cdots c_{\lm_1}
\underbrace{v_1\otimes\cdots\otimes v_1}_{\lm_1\mbox{ times }}
\otimes\cdots\otimes
\underbrace{v_n\otimes\cdots\otimes v_n}_{\lm_n\mbox{ times }},
$$
where $c_i$ is the sum of the elements $\sgn(\sg)\sg$
taken over $\sg\in G$ such that $\sg(i)=i$ for any $i$
not belonging to
$\{i+\sum_{j=1}^s\lm_j:0\le s\le n-1,\, i\le\lm_{s+1}\}$.
It is obvious that $Av(\lm)=KGv(\lm)$ is isomorphic as a $KG$-module to the
Specht module $S^\lm$ defined in~\cite{James1}.
Indeed, let us identify an element
$v_{i_1}\otimes$\cc$\otimes v_{i_r}$ of $(E^{\otimes r})^\lm$
with the tabloid that has each $s$, $s=1$,\ldots,$r$
in row $i_s$.
This correspondence extended linearly to
$(E^{\otimes r})^\lm$ gives an isomorphism of $KG$-modules
$(E^{\otimes r})^\lm$ and $M^\lm$.
Tabloids and $M^\lm$ are understood in the sense
of~\cite[Definitions~3.9 and~4.1]{James1}
with the only difference that in our case the symmetric group acts
on the left.
Clearly, the restriction of this isomorphism to $Av(\lm)$ gives the required
correspondence.

Let $\lm$ be an arbitrary partition of $X$.
Denote by $P$ the maximal $A$-submodule of $Av(t(\lm))$.
The modules $Av(t(\lm))/P$ and $D^{t(\lm)}$ are isomorphic as
$KG$-modules and therefore are isomorphic as $A$-modules.
It follows from this fact and from $e_\lm D^{t(\lm)}\ne0$
that $e_\lm Av(t(\lm))\ne0$ and $e_iAv(t(\lm))\ne0$ for some $i\in I_\lm$.
All elements of the nonzero space $e_iAv(t(\lm))$ have weight $t(\lm)$.
On the other hand $e_iAv(t(\lm))\st M_i\cong P(\lm)$.
Hence $t(\lm)$ is a weight of $P(\lm)$ and $t(\lm)\tleq\lm$
by property~\ref{p:3}.

Let us prove the reverse inequality.
Without loss of generality we may assume
$I_\lm=\{1$,\ldots,$l\}$.
Here $l\ge1$ since $e_\lm\ne0$.
As is well known, $KGL_n(K)v(\lm)$ is isomorphic to
the so-called Weyl module with highest weight $\lm$.
But this same module is embeddable in $P(\lm)$.
Thus there is an embedding of $GL_n(K)$-modules
$\tau:KGL_n(K)v(\lm)\to M_1$.
Let $v=\tau(v(\lm))$.

Let $U$ be the subgroup of $GL_n(K)$ generated by
endomorphisms $\vp$ such that $\vp(v_k)=v_k$ for $k\ne i$ and
$\vp(v_i)=v_i+\alpha v_j$, where $\alpha\in K$ and $j<i$.
Corollary~17.18 of~\cite{James1} actually means that an element
$m\in(E^{\otimes r})^\lm$ belongs to $Av(\lm)$ if and only if
$um=m$ for any $u\in U$.
Using this fact and that $\tau$ is an embedding of $GL_n(K)$-modules,
we get $v\in Av(\lm)$.

If $e_\lm$ had annihilated each composition factor of
$Av$, we would have got $e_\lm Av=(e_\lm)^sAv=0$.
However $v\in e_\lm Av$.
Hence $e_\lm D^\mu\ne0$, where
$D^\mu$ is a composition factor of $Av$
and therefore a composition factor of $Av(\lm)$.
Since $Av(\lm)\cong S^\lm$, we have $\mu\treq \lm$.
By definition $\mu=t(\lm)$ and therefore $\mu\tleq\lm$
by what was proved above.
Hence we obtain $t(\lm)=\lm$ and $Av=Av(\lm)$.
\endproof

During the proof we additionally obtained
\begin{sled}
Let $v$ be a nonzero element of $M_i\cong P(\lm)$ of weight $\lm$.
Then the $GL_n(K)$-module generated by $v$ is isomorphic to the Weyl module
with highest weight $\lm$ and the $G(r)$-module
generated by $v$ is isomorphic to the Specht module $S^\lm$.
\end{sled}

Lemma~\ref{l:0.25} is very important for our theory.
It gives a description via the primitive idempotents of $A$ of
the same irreducible modules that were described in~\cite{James1}
and that were studied in~\cite{Shchigolev9}.
Without this lemma the approach taken here would be ineffective.
The author does not know for sure how to cite this assertion.
Therefore its proof is given in this article.

\section{Embedding operators}\label{s:2}

In the previous section we chose an arbitrary decomposition of
$E^{\otimes r}$ into a direct sum of irreducible tilting modules.
However, it is convenient that these decompositions for different $r$
be consistent.
Thus we will build them inductively as follows.

First we put $M_1(1)=E$.
Suppose that we have a decomposition
$E^{\otimes r}=\bigoplus_{i=1}^{h_r}M_i(r)$.
By properties~\ref{p:1} and~\ref{p:2} each module $M_i(r)\otimes E$
is tilting.
Therefore by property~\ref{p:4} it is decomposable into a
direct sum of indecomposable tilting modules.
Collecting all these summands, we
obtain a decomposition of $E^{\otimes r+1}$.

Let $r<r'$ be two positive integers.
Let us define an embedding operator $f_{r',r}:A(r) \to A(r')$
by the formula $f_{r',r}(\vp)=\vp \otimes i_{E^{\otimes r'-r}}$.
Then we have $f_{r',r}\circ\sg_r=\sg_{r'}|_{KG(r)}$ and
$f_{r'',r'}\circ f_{r',r}=f_{r'',r}$.
It is obvious that
$M_i(r)\otimes E^{\otimes r'-r}=\bigoplus_{j\in X}M_j(r')$
implies $f_{r',r}(e_i^r)=\sum_{j\in X}e_j^{r'}$.

Let us recall the terminology and the definitions used in~\cite{Shchigolev9}.
A partition $\lm\in\Lm^+_r$ is called {\it degenerate}
({\it $n$-degenerate}) if $\lm_n=0$.
For two partitions $\lm,\mu\in\Lm^+_r$ and
two integers $\alpha$, $\beta$, we denote by $\alpha\lm+\beta\mu$
the componentwise linear combination of these partitions.
Let $(i^n)$ denote the sequence of length $n$, whose every entry is $i$.
Denote by $d_r$ the sum of the idempotents $e_\lm$ over all
degenerate $\lm\in\Lm^+_r$.

Now suppose that $\lm$ is a partition of $\Lm^+_r$ satisfying
\begin{usl}\label{usl:1}
There exists $a$ such that
for any degenerate $\mu\in\Lm^+_{r+an}$ and $k\ge a$
the module $P(\mu)\otimes E^{\otimes(k-a)n}$ does not have
a direct summand isomorphic to $P(\lm+(k^n))$.
\end{usl}
The absence of the direct summand mentioned above
for any degenerate $\mu$ is equivalent to
$e_{\lm+(k^n)}f_{r+kn,r+an}(d_{r+an})=0$.
We will see in Section~\ref{s:3} that this assumption is not meaningless.

\begin{teo}\label{t:1}
Suppose that $\lm+(1^n)$ is a $p$-regular partition.
Let $k\ge \frac{r^2}n+(2r+1)a+a^2n$
and $\mu$ be a degenerate $p$-regular partition of $\Lm^+_{r+kn}$.
Then any $KG(r+kn)$-module $M$ of length $2$
with factors $D^{\lm+(k^n)}$ and $D^\mu$ is decomposable.
\end{teo}
{\bf Proof.}
Considering, if necessary, the dual modules instead of the initial ones,
we may suppose without loss of generality that there exists a
submodule $N\st M$ such that $N\cong D^\mu$ and
$M/N\cong D^{\lm+(k^n)}$.

Let $x$ and $y$ be elements of $KG(r+an)$ such that
$\sg_{r+an}(x)=1-d_{r+an}$ and $\sg_{r+an}(y)=d_{r+an}$.
In addition we have $\sg_{r+kn}(x)=f_{r+kn,r+an}(1-d_{r+an})$ and
$\sg_{r+kn}(y)=f_{r+kn,r+an}(d_{r+an})$.

We have the following equalities:
\begin{multline}\label{eq:1}
\shoveleft{%
\sg_{r+kn}(y)e_{\lm+(k^n)}=0,\,\sg_{r+kn}(x)e_\mu=0,\,
\sg_{r+kn}(y)e_\mu=e_\mu,\,\sg_{r+kn}(x)e_{\lm+(k^n)}=e_{\lm+(k^n)}.}
\end{multline}
The first one follows from the assumption about $\lm$ we have made.
The second one follows from the fact that for any nondegenerate
$\nu\in\Lm^+_{r+an}$ the module $P(\nu)\otimes E^{\otimes(k-a)n}$
does not have a direct summand isomorphic to $P(\mu)$.
The last two equalities follow from the two equalities already proved and
from the formula $\sg_{r+kn}(x)+\sg_{r+kn}(y)=1$.

Consider the translates of $x$ defined by the formula
$x_i=\pi_ix\pi_i^{-1}$, $i=1$, \ldots ,$r+an+1$,
where $\pi_i=\prod_{j=1}^{r+an}(j+(i-1)(r+an),j)$.
Note that $\supp(x_i)=\{1+(i-1)(r+an)$, \ldots, $i(r+an)\}$.
Since $\supp(x_i)$ are disjoint for distinct $i$,
the elements $x_i$ commute.
The lower bound for $k$ in the hypothesis of the theorem
makes all $x_i$ belong to $KG(r+kn)$.

We have the relations
\begin{equation}\label{eq:2}
 yM\st N, \quad x_iN=0, \quad yN\ne0,
\end{equation}
which follow from the first two equalities~(\ref{eq:1})
and Lemma~\ref{l:0.5}.

Choose an arbitrary $m\in M\setminus N$ and prove that
$x_1$\cc$x_{r+an+1}m\notin N$.
Assume the contrary.
Then $\sg_{r+kn}(x_1)$\cc$\sg_{r+kn}(x_{r+an+1})\bar m=0$,
where $\bar m$ is the image of $m$ under the projection $M\to M/N$.
Hence by the last formula of~(\ref{eq:1}) and Lemma~\ref{l:0.5},
we obtain
$$
\begin{array}{l}
 \sg_{r+kn}(x_i)\bar m=
 \sg_{r+kn}(\pi_i)\sg_{r+kn}(x)\sg_{r+kn}%
 (\pi_i^{-1})\bar m=\\[12pt]
 \sg_{r+kn}(\pi_i)\sg_{r+kn}(x)e_{\lm+(k^n)}%
 \sg_{r+kn}(\pi_i^{-1})\bar m=\bar m.
\end{array}
$$
This implies $\bar m=0$, which leads to a contradiction.

Denote by $P$ the $KG(r+kn)$-module generated by
$x_1$\cc$x_{r+an+1}m$.
Let us calculate what $yP$ is.
Let $\alpha\in G(r+kn)$.
The cardinality of $\supp(\alpha^{-1}y\alpha)$ does not exceed $r+an$.
Thus there is an $i=1$, \ldots, $r+an+1$ such that
$\alpha^{-1}y\alpha$ and $x_i$ commute.
We have
$$
y\alpha x_1\cdots x_{r+an+1}m\in\alpha x_i(\alpha^{-1}y\alpha)M\st
\alpha x_iN=0.
$$
Thus $yP=0$.
Since $yN\ne0$, we have $P\cap N=0$ and $M=P\oplus N$.
\endproof

By induction on the length one can easily prove that under
the hypothesis of Theorem~\ref{t:1}
any $KG(r+kn)$-module $M$ of a finite length such that one
of its composition factors is isomorphic to $D^{\lm+(k^n)}$ and
the remaining factors are degenerate, contains a submodule
isomorphic to $D^{\lm+(k^n)}$.

\section{Examples}\label{s:3}

We proved Theorem~\ref{t:1} under the assumption that $\lm$
satisfies Condition~\ref{usl:1}.
In this section we are going to give two examples of such partitions.

Example~1.
Suppose $n<p$.
For any positive integer $R$ we put
$C_0(R)=\{\lm\in\Lm^+_R:\lm_1-\lm_n\le p-n\}$.
Lemma~12(2) of~\cite{Mathieu1} asserts that in the case where
$\mu\in\Lm_R^+\setminus C_0(R)$, the module
$P(\mu)\otimes E$ does not have a direct summand isomorphic to
$P(\nu)$, where $\nu\in C_0(R+1)$.
It obviously follows from this fact that for any $i$
the module $P(\mu)\otimes E^i$ does not have a direct summand
isomorphic to $P(\nu)$, where $\nu\in C_0(R+i)$.

Choose an arbitrary $\lm$ of $C_0(r)$ and
prove that this partition satisfies Condition~\ref{usl:1}.

\begin{lemma}\label{l:1}
Let $R\ge(n-1)(p-n)+1$ and $\mu$ be a degenerate partition of $\Lm_R^+$.
Then $\mu\notin C_0(R)$.
\end{lemma}
{\bf Proof.}
We have $(n-1)\mu_1\ge\mu_1+$\cc$+\mu_{n-1}=R>(n-1)(p-n)$.
Hence $\mu_1-\mu_n=\mu_1>p-n$ and $\mu\notin C_0(R)$.
\endproof

Choose any, for example minimal, $a$ satisfying the inequality
$r+an\ge(n-1)(p-n)+1$.
Then by~\cite[Lemma~12(2)]{Mathieu1} and Lemma~\ref{l:1} the partition
$\lm$ satisfies Condition~\ref{usl:1} for $a$ defined above.
Thus Theorem~\ref{t:1} is applicable to $\lm$.
However taking into account that we consider the special case $n<p$,
the estimate of this theorem can be sharpened.
\begin{teo}\label{t:2}
Let $\lm\in C_0(r)$, $\lm_n\ge(n-1)(p-n)+2$
and $\mu$ be a degenerate partition of $\Lm^+_r$.
Then any $KG(r)$-module $M$ of length $2$
with factors $D^\lm$ and $D^\mu$ is decomposable.
\end{teo}
{\bf Proof}
differs from that of Theorem~\ref{t:1} mainly
by the choice of $x$ and $y$.
Similarly to Theorem~\ref{t:1} we assume without loss of generality that
there is a submodule $N\st M$ such that $N\cong D^\mu$ and
$M/N\cong D^\lm$.
For brevity, we put $R=(n-1)(p-n)+1$.

Let $x=\frac1{n!}\sum_{\sg\in G(n)}\sgn(\sg)\sg$ and
$x_i=\pi_ig\pi_i^{-1}$, $i=1$, \ldots ,$R+1$,
where $\pi_i=\prod_{j=1}^n (j+(i-1)n,j)$.
Choose $y\in KG(R)$ so as to have $\sg_R(y)=d_R$.
Then $\sg_r(y)=f_{r,R}(d_R)$.
By~\cite[Lemma~12(2)]{Mathieu1} and Lemma~\ref{l:1}
we obtain $\sg_r(y)e_\lm=0$.
The first and the third of formulas~(\ref{eq:2}) can be checked
similarly to what was done in Theorem~\ref{t:1}.
The second formula follows from the fact that $\mu$
contains no more than $n-1$ nonzero parts.
From the direct construction of $D^\lm$ as a quotient module of
the Specht module $S^\lm$ and from the inequality $\lm_n\ge R+1$ it follows
that there exists some $m\in M$ such that $x_1$\cc$x_{R+1}m\notin N$.
Denote by $P$ the $KG(r)$-submodule of $M$ generated by
$x_1$\cc$x_{R+1}m$ and similarly to Theorem~\ref{t:1} show that
$yP=0$.
Hence $M=P\oplus N$.
\endproof

Example~2.
Suppose $n=2$.
Choose an arbitrary $\lm\in\Lm^+_r$.
Let $\mu\in\Lm^+_{r+1}$ be a partition such that
$P(\mu)$ is a direct summand of $P(\lm)\otimes E$.
Suppose that $k$ is a positive integer such that $\lm_1-\lm_2\ge p^k-1$.
Let us prove that in this case $\mu_1-\mu_2\ge p^k-1$.
The proof proceeds by applying the formulas for the product
of tilting modules form~\cite[Section~1]{Erdmann1}.
These formulas deal with tilting modules over $SL_2(K)$.
However, it is clear that if a true formula of the form
$T(s)\otimes T(1)=T(s')$
is replaced by
$P((\frac{r+s}2,\frac{r-s}2))\otimes E=P((\frac{r+1+s'}2,\frac{r+1-s'}2))$,
then we have a true formula again.
In what follows we will take account of such reformulations
while citing formulas of~\cite{Erdmann1}.

For brevity, we put $m=\lm_1-\lm_2$ and $m'=\mu_1-\mu_2$.
We assume $m>0$ and $k>0$, since in the contrary case
the assertion being proved is obvious.
Let $m=\sum_{s=0}^\infty i_sp^s$, where $0\le i_s\le p-1$,
be the $p$-addic expansion of $m$.

If $i_0=p-1$ then by~\cite[Lemma~1.5(1)]{Erdmann1} we have
$m'=m+1\ge p^k$.
If $0\le i_0\le p-3$ (the case $p>2$) then $m\ge p^k$.
By~\cite[Lemma~1.5(2),(3)]{Erdmann1} we have
$m'=m+1$ or $m'=m-1$.
In both cases we have $m'\ge p^k-1$.

Now suppose that $i_0=p-2$.
Denote by $t$ the positive integer such that
$i_s=p-1$ for $s=1$, \ldots, $t-1$ and $i_t<p-1$.
We put $\sg=\sum_{s=t+1}^\infty i_sp^{s-t-1}$.
Then $m+2=(i_t+1)p^t+\sg p^{t+1}$.

Case $p>2$.
For the calculation we will use \cite[Lemma~1.7.3]{Erdmann1}.
To prove the inequality $m'\ge p^k-1$
it suffices to prove that $m+1-2p^{t-1}\ge p^k-1$
and that $m+1-2p^t\ge p^k-1$ for $i_t=1$ and $\sg>0$
for $i_t\ge2$.

Suppose that $t<k$.
Since $m\ge p^k$, we have $i_{k'}>0$ for some $k'\ge k$.
Hence $m+2\ge(i_t+1)p^t+p^k$ and $m+1-2p^{t-1}\ge p^k-1+p^t-2p^{t-1}>p^k-1$.
If $i_t\ge1$ then $m+1-2p^t\ge p^k-1+(i_t-1)p^t\ge p^k-1$.

Suppose that $t=k$.
Since $m\ge p^k$, we have $i_k>0$ and $\sg=0$ or $\sg>0$.
In the first case $m+2\ge2p^k$ and $m+1-2p^{t-1}\ge 2p^k-1-2p^{k-1}>p^k-1$.
In addition, if $i_k\ge2$ then $m+2\ge3p^k$ and $m+1-2p^t\ge p^k-1$.
In the second case $m+2\ge3p^k$ and $m+1-2p^t\ge p^k-1$.

Suppose $t>k$.
Since $m+2\ge p^t$, we have
 $m+1-2p^{t-1}\ge p^{t-1}-1\ge p^k-1$.
If $\sg>0$ or $i_t\ge2$, then $m+2\ge 3p^t$ and
$m+1-2p^t\ge p^t-1>p^k-1$.

Case $p=2$.
For the calculation we will use~\cite[Lemma~1.7.2]{Erdmann1}.
In this case to prove the inequality $m'\ge 2^k-1$
it suffices to prove that $m+1-2^{t-1}\ge 2^k-1$
and that $m+1-2^t\ge 2^k-1$ for $\sg>0$.

Suppose $t<k$.
Similarly to the case $p>2$ we obtain $m+2\ge2^t+2^k$.
Hence $m+1-2^t\ge 2^k-1$.

Suppose $t=k$.
Similarly to the case $p>2$ we obtain $i_k>0$ or $\sg>0$.
In both cases we have $m+2\ge2^{k+1}$ and $m+1-2^t\ge 2^{k+1}-1-2^k=2^k-1$.

Suppose $t>k$.
Since $m+2\ge2^t$, we have $m+1-2^{t-1}\ge 2^{t-1}-1\ge 2^k-1$.
If $\sg>0$ then $m+2\ge 2^{t+1}$ and
$m+1-2^t\ge 2^t-1>2^k-1$.

Now we are ready to check Condition~\ref{usl:1} for
any partition $\lm\in\Lm^+_r$.
Choose any, for example minimal, $k$
such that $\lm_1-\lm_2<p^k-1$ and put $a=p^k-1-r$.
From the above calculations one can see that Condition~\ref{usl:1}
is satisfied for $a$ we have chosen.

\section{Finite basis property}\label{s:4}

Let us study what consequences of the results of Section~\ref{s:3}
one can get as far as the operation $\upa$ introduced in~\cite{Shchigolev9}
is concerned.

For any $KG(r)$-module $V$ we denote by $\dual(V)$ the dual module
of $V$.
In addition, there exists the nondegenerate form
$\<\:,\,\>:V\times \dual(V)\to K$ such that
$\<\sg v,u\>=\<v,\sg^{-1} u\>$ for any
$v\in V$, $u\in\dual(V)$ and $\sg\in G(r)$.
For a subset $A$ of $V$ or of $\dual(V)$ we denote by $V^\perp$
the set $\{u\in\dual(V):\<A,u\>=0\}$ or the set $\{v\in V:\<v,A\>=0\}$
respectively.

\begin{lemma}\label{l:2}
Let $\lm$ be a $p$-regular nondegenerate partition of $\Lm^+_r$
and $V$ be a proper submodule of the Specht module $S^\lm$.
Consider the following conditions:

\begin{enumerate}
\item\label{pp:1} $V\doa\upa\stn V$.
\item\label{pp:2}
There is a $KG(r)$-module $L$ and a nonisomorphic embedding
{$\iota:\dual(S^\lm/V)\to L$} such that
$\im\iota$ is an essential submodule of $L$
and all composition factors of $L/\im\iota$ are degenerate.
\end{enumerate}

Condition~\ref{pp:2} always follows from Condition~\ref{pp:1}.
The reverse implication is true for $n<p$.
\end{lemma}
{\bf Proof.}
\ref{pp:1}$\Rightarrow$\ref{pp:2}
Define $L=\dual(S^\lm/V\doa\upa)$ and $\iota$
to be the embedding induced by the projection
$S^\lm/V\doa\upa\to S^\lm/V$.
From the obvious equality $V\doa\upa\doa=V\doa$
(see Introduction) and~\cite[Theorem~5(a)]{Shchigolev9} it follows that
$V/V\doa\upa$.
Therefore all composition factors of $L/\im\iota$ are also degenerate.
It remains to show that $\im\iota$ is an essential submodule of $L$.
Let $N$ be a submodule of $L$ and $N\cap\im\iota=0$.
Then we have $N^\perp+(V/V\doa\upa)=S^\lm/V\doa\upa$.
Since $P^\lm/V\doa\upa$ is a unique maximal submodule of
$S^\lm/V\doa\upa$ and $V\st P^\lm$, we have $N^\perp=S^\lm/V\doa\upa$.
Hence $N=0$.

\ref{pp:2}$\Rightarrow$\ref{pp:1} subject to $n<p$.
Without loss of generality we may assume that
$L/\im\iota\cong D^\mu$, where $\mu$ is a degenerate
$p$-regular partition.

In the remaining part of the proof we will follow the definitions
of~\cite{J.A.Green1}.
Choose an arbitrary $N\ge r$ and put $G=G(r)$, $S=S_K(N,r)$,
$\om=(1^r,0^{N-r})$, $e=\xi_\om$.
Let us fix the sequence $u(r)=(1$, \ldots, $r)$.
The map $\xi_{u(r)\pi,u(r)}\to\pi$ defines an isomorphisms of
the rings $eSe$ and $KG$.
Thus we may consider any $KG$-module as an $eSe$-module and vice versa.

Let $V_{\lm',K}$ be the Weyl module with highest weight $\lm'$.
In the paper~\cite{J.A.Green1} the isomorphism~(6.3d)
was used to establish an isomorphism of $V_{\lm',K}^\om$
and $\dual(S^{\lm'})$
($\bar S_{T,K}$ in the notation of~\cite{J.A.Green1}).

Let $K_s$ denote the sign representation of $G$.
There is the standard embedding of $V\otimes K_s$ into
$\dual(S^{\lm'})$ and therefore into $V_{\lm',K}^\om$.
Denote by $V'$ the image of $V$ under this embedding.

Put $U=SV'$.
This module is a submodule of $V_{\lm',K}$.
The paper~\cite{J.A.Green1} introduces the left
$S$-modules $D_{\lm',K}$ and $A_K(N,r)$.
The last module is also a right $S$-module.
In addition $D_{\lm',K}$ is a submodule of $\^{\lm'}A_K(N,r)$
that is a right weight space of $A_K(N,r)$.

According to~\cite[Section~5.1]{J.A.Green1}
the restriction of the natural form $\<\:,\,\>$
defined on $E^{\otimes r}$ gives a nondegenerate contravariant form
$(\,,\,):V_{\lm',K}\times D_{\lm',K}\to K$.
Put $W=\{w\in D_{\lm',K}:(U,w)=0\}$.
The module $W$ is the contravariant dual to $V_{\lm',K}/U$
in the sense of~\cite[Section~2.7]{J.A.Green1}.
Therefore $eW$ is dual to $V_{\lm',K}^\om/V'$
as a $KG$-module in the usual sense (see~\cite[page 92 ]{J.A.Green1}).
In addition $V_{\lm',K}^\om/V'\cong S^\lm/V\otimes K_s$.
Thus $eW\cong\dual(S^\lm/V)\otimes K_s$.
Consider any isomorphism $\theta:eW\to\dual(S^\lm/V)\otimes K_s$.
We put $M=L\otimes K_s$ and $\vp=(\iota\otimes i_{K_s})\theta$.
Clearly, $\vp$ is an embedding of $eW$ to $M$.
In addition, $\im\vp$ is an essential submodule of $M$
and $M/\im\vp\cong D^\mu\otimes K_s$.

Put $h(M)=Se\otimes_{eSe}M$.
Denote by $W'$ the $K$-subspace of $h(M)$
spanned by the elements of the form $se\otimes\vp(ew)$,
where $s\in S$ and $w\in W$.
Clearly, $W'$ is an $S$-module.
Consider the map $\psi:W'\to W$ defined by the formula
$\psi(se\otimes\vp(ew))=sew$.

Since $W\st\^{\lm'}A_K(N,r)$ and the latter module is injective,
there exists an extension
$\bar\psi:h(M)\to\^{\lm'}A_K(N,r)$ of $\psi$.
Consider the isomorphism $\ep:M\to eh(M)$
defined by the formula $\ep(m)=e\otimes m$.
Let $\tau=\bar\psi\ep$ be their composition.
It is obvious that $\tau\vp=i_{eW}$.
Therefore $\tau$ isomorphically takes $\im\vp$ to $eW$.
Since $\im\vp$ is an essential submodule of $M$,
we have that $\tau$ is an embedding of $M$ in $e\,\^{\lm'}A_K(N,r)$.

The module $e\,\^{\lm'}A_K(N,r)$ can be identified with
$M^{\lm'}$ in the notation of~\cite{James1}.
Then $eD_{\lm',K}$ is identified with the Specht module $S^{\lm'}$.
Applying the results of~\cite[Section~17]{James1}, we obtain
that there is a filtration of $e\,\^{\lm'}A_K(N,r)/eD_{\lm',K}$
with factors isomorphic to $S^\nu$, where $\nu\triangleright\lm'$.
Therefore either $\tau(M)\st eD_{\lm',K}$ or
some $S^\nu$, where $\nu\vartriangleright\lm'$,
contains a submodule isomorphic to $D^{\mu}\otimes K_s$.

Let us show that the latter case is impossible.
Since $S^\nu\st M^\nu$ (see~\cite{James1})
and the modules $D^{\mu}\otimes K_s$ and $M^\nu$ are self-dual,
there exists an epimorphism $\pi:M^\nu\to D^{\mu}\otimes K_s$.
Put $g=\sum_{\sg\in G(n)}\sg$.
Since $\nu\vartriangleright\lm'$, we have $\nu_1\ge n$.
It follows from this fact and from $n<p$ that $M^\nu$ has
a cyclic generator of the form $gv$, where $v\in M^\nu$.
Applying $\pi$, we obtain $\pi(gv)=g\pi(v)\in g(D^{\mu}\otimes K_s)=0$.
The last equality follows from the fact that $\mu$ is degenerate.
This leads to a contradiction.

Thus we have obtained that $eW\st\tau(M) \st eD_{\lm',K}$.
Applying the form $(\:,\,)$ and multiplying by $K_s$,
we obtain that there exists some submodule $V_0$
of $V$ such that $V/V_0\cong D^\mu$.
By~\cite[Corollary~6]{Shchigolev9}we have $V\doa\upa\stn V$.
\endproof

{\bf Note.}
In the first part of the proof it is possible to consider an
arbitrary field, while in the second part the field should be infinite.
However, the second part is given only for an illustration
of the relation between properties~\ref{pp:1} and~\ref{pp:2}
and is never used in the present paper.

For any $p$-regular partition $\nu$, we denote by $P^\nu$
the unique maximal submodule of the Specht module $S^\nu$.
One can easily check that if $\nu$ is nondegenerate, then
$P^\nu\doa=P^{\nu-(1^n)}$.
Now by the first part of Lemma~\ref{l:2}, we have
that under the hypothesis of Theorems~\ref{t:1} and~\ref{t:2}
the equalities
\begin{equation}\label{eq:3}
P^{\lm+(k-1^n)}\upa=P^{\lm+(k^n)}, \quad P^{\lm-(1^n)}\upa=P^\lm
\end{equation}
hold respectively.
If $K$ is not algebraically close, then applying the natural isomorphism
$S_K^\nu\otimes\bar K\cong S_{\bar K}^\nu$
which takes $P_K^\nu\otimes\bar K$ to $P_{\bar K}^\nu$
we obtain formulas~(\ref{eq:3}) for any field of the same characteristic.
It follows from these formulas that for any $\lm$ satisfying
Condition~\ref{usl:1}, there exists finitely many relations
defining all irreducible modules
$D^{\lm+(k^n)}\cong S^{\lm+(k^n)}/P^{\lm+(k^n)}$,
$k\in\N$.

\bibliography{ref}
\bibliographystyle{m}

\end{document}